\documentclass[11pt]{amsart}

\input xy
\xyoption{all}

\date{August 30, 2004}

\begin{document}

\title{Invariant chains and the homology of quotient spaces}

\author{Kevin P.~Knudson}\thanks{Partially supported by NSF grant no.~DMS-0242906 and by ORAU}

\address{Department of Mathematics and Statistics, Mississippi State University, Mississippi State, MS 39762}
\email{knudson@math.msstate.edu}

\newtheorem{theorem}{Theorem}[section]
\newtheorem{prop}[theorem]{Proposition}
\newtheorem{lemma}[theorem]{Lemma}
\newtheorem{cor}[theorem]{Corollary}
\newtheorem{conj}[theorem]{Conjecture}
\newtheorem{definition}[theorem]{Definition}
\newtheorem{remark}[theorem]{Remark}

\newcommand{\zz}{{\mathbb Z}}
\newcommand{\zq}{{\mathbb Q}}
\newcommand{\ra}{\rightarrow}
\newcommand{\lra}{\longrightarrow}
\newcommand{\bop}{\bigoplus}
\newcommand{\cp}{{\mathbb C}{\mathbb P}}
\newcommand{\rp}{{\mathbb R}{\mathbb P}}
\newcommand{\zr}{{\mathbb R}}
\newcommand{\zc}{{\mathbb C}}

\begin{abstract}  For a finite group $G$ and a finite $G$-CW-complex $X$, we construct groups
$H_\bullet(G,X)$ as the homology groups of the $G$-invariants of
the cellular chain complex $C_\bullet(X)$. These groups are
related to the homology of the quotient space $X/G$ via a norm
map, and therefore provide a mechanism for calculating
$H_\bullet(X/G)$.  We compute several examples and provide a new
proof of ``Smith theory": if $G=\zz/p$ and $X$ is a mod $p$
homology sphere on which $G$ acts, then the subcomplex $X^G$ is
empty or a mod $p$ homology sphere.  We also get a new proof of
the Conner conjecture:  If $G=\zz/p$ acts on a $\zz$-acyclic space
$X$, then $X/G$ is $\zz$-acyclic.
\end{abstract}

\maketitle

\section{Introduction}

There is a rich literature devoted to the problem of computing the (co)-homology of quotient spaces.  Indeed,
this question is but a small part of the field of equivariant topology.  Given a finite group $G$ and a Mackey
functor $M$, there is a cohomology theory $H^\bullet_G(-;M)$ (ordinary Bredon cohomology) with the property
that if $X$ is a $G$-CW-complex, then
$$H^\bullet_G(X;\underline{\zz}) \cong H^\bullet(X/G;\zz),$$ where $\underline{\zz}$ is the constant Mackey functor
associated to $\zz$.  The usual universal coefficient theorems hold, and there are long exact sequences associated
to short exact sequences of Mackey functors.  In the special case of $G=\zz/p$, one can derive useful exact
sequences relating $H^\bullet_G(X;\zz/p)$ with $H^\bullet(X^G;\zz/p)$ and $H^\bullet(X-X^G;\zz/p)$ \cite{borel}.  What's
more, Bredon cohomology can be generalized to the case where $G$ is a compact Lie group.

In this paper, we take a different approach to the homology of the quotient space $X/G$.  The homology of this
space is computed via the chain complex $C_\bullet(X)_G$ of coinvariants of the $G$-action on the cellular chain complex
of $X$.  The norm map $$N:C_\bullet(X)_G\lra C_\bullet(X)^G$$ defined by
$$N(\overline{\sigma})=\sum_{g\in G}g\sigma$$ ($\overline{\sigma}$ denotes the $G$-orbit of the cell $\sigma$) is
an injective map of chain complexes.  Note that $N$ is {\em not} an isomorphism, although it becomes one upon inverting $|G|$.
Denote the homology groups of $C_\bullet(X)^G$ by $H_\bullet(G,X)$.  In the case $G=\zz/p$, these groups are related to the
groups $H_\bullet(X/G)$ as follows.

\medskip

\noindent {\bf Proposition \ref{lesprop}.} {\em If $G=\zz/p$, then there is a long exact sequence
$$\ra H_n(X/G;\zz)\ra H_n(G,X;\zz)\ra H_n(X^G;\zz/p)\ra H_{n-1}(X/G;\zz)\ra.$$}

In simple examples, the groups $H_\bullet(G,X)$ are easy to compute, but then so are the $H_\bullet(X/G)$.  What
one needs is another mechanism to compute the $H_\bullet(G,X)$.  To obtain this, we proceed as follows.  Note
that the invariants functor $(-)^G$ is the $0$-th cohomology of $G$.  Let $P_\bullet\ra\zz$ be a projective resolution
over $\zz G$ and consider the fourth quadrant double complex
$$E^0_{s,t} = \text{Hom}_{\zz G}(P_{-t},C_s(X)).$$ Define $S_n(G,X)$ to be the $n$-th homology of the total complex
of $E^0_{\bullet,\bullet}$.  As usual, one obtains two spectral sequences converging to $S_\bullet(G,X)$.  One of these
has $$E^2_{s,0}=H_s(G,X)$$ and so one may glean information about $H_\bullet(G,X)$ via these auxiliary functors $S_\bullet(G,X)$.

As an application of these techniques, we obtain a new proof of ``Smith theory".

\medskip

\noindent {\bf Theorem \ref{smiththm}.} {\em Let $X$ be a finite-dimensional CW-complex on which $G=\zz/p$ acts
cellularly.  Then

\begin{enumerate}
\item If $H_\bullet(X;\zz/p)\cong H_\bullet(S^n;\zz/p)$ for some $n$, then either $X^G=\emptyset$ or
$H_\bullet(X^G;\zz/p)\cong H_\bullet(S^m;\zz/p)$ for some $m$; and

\item If $H_\bullet(X;\zz/p)\cong H_\bullet(\text{\em pt};\zz/p)$, then $H_\bullet(X^G;\zz/p)\cong
H_\bullet(\text{\em pt};\zz/p)$.

\end{enumerate}}

\medskip

We also get a proof of the Conner conjecture (see \cite{oliver}).

\medskip

\noindent {\bf Theorem \ref{conner}.}  {\em Let $G=\zz/p$ and let
$X$ a finite $G$-CW-complex.  If $X$ is $\zz$-acyclic, then so is
$X/G$.}

\medskip

In addition, we compute some examples in Section \ref{examples}.
For instance, we obtain a calculation of
$H_\bullet(\cp^\infty/(\zz/2);\zz)$ where $\zz/2$ acts via complex
conjugation.

Although we do not compute $H_\bullet(X/G)$ for any new examples,
it is our hope that the techniques developed here will prove
useful for calculations. We should also point out that the
approach taken here is similar in spirit to a construction of
S.~Waner \cite{waner}, who showed that Bredon cohomology can be
computed by using the Mackey functor $G/H\mapsto C_\bullet(X)^H$.
Our method is quite different, however.

\medskip

\noindent {\em Acknowledgements.}  I thank Fred Cohen for many useful conversations, and in particular for pointing
out reference \cite{mitchell} to me.

\section{Basic Definitions}\label{basic}

In this paper, $G$ is a finite group and $X$ is a $G$-CW-complex.  This means that $X$ is a CW-complex equipped with
a cellular $G$-action such that if $g\in G$ fixes a cell $\sigma$, then $g$ fixes $\sigma$ pointwise.  It follows that
the fixed point set $X^G$ is a subcomplex and the quotient space $X/G$ is a CW-complex.

Denote by $C_\bullet(X;A)$ the cellular chain complex of $X$ with coefficients in the abelian group $A$.  A typical element
of $C_i(X;A)$ is a linear combination of elements of the form $\sigma\otimes a$ with $\sigma$ an $i$-cell of $X$ and $a\in A$.
If $A=\zz$, we usually drop it from the notation.

The action of $G$ on $X$ induces an action on $C_\bullet(X;A)$ by $$g(\sigma\otimes a)\mapsto (g\sigma)\otimes a.$$ Denote
the subcomplex of $G$-invariants by $C_\bullet(X;A)^G$.  Note that $C_\bullet(X^G;A)$ is a {\em proper} subcomplex
of $C_\bullet(X;A)^G$.

\begin{definition} For all $i\ge 0$, define $$H_i(G,X;A)=h_i(C_\bullet(X;A)^G).$$
\end{definition}

\subsection{Variance}\label{variance} The functors $H_i(G,X;A)$ are covariant in $X$ and contravariant in $G$; the notation
was chosen to reflect this (compare with Hom).  Precisely, if $f:X\ra Y$ is a cellular map of $G$-CW-complexes, then there
is a map
$$f_*:H_i(G,X;A)\lra H_i(G,Y;A)$$ defined as follows.  Consider the chain map $f_*:C_i(X;A)\ra C_i(Y;A)$ given by
$\sigma\otimes a\mapsto f(\sigma)\otimes a$.  If $z\in C_i(X;A)^G$, then if $g\in G$,
$$gf_*(z) = f_*(gz)=f_*(z)$$ so that $f_*$ restricts to a chain map
$$f_*:C_\bullet(X;A)^G\lra C_\bullet(Y;A)^G.$$

If $\alpha:G'\ra G$ is a homomorphism and $X$ is a $G$-CW-complex, then $X$ is a $G'$-complex via $\alpha$.  In this case, there
is a homomorphism $$\alpha^*:H_i(G,X;A)\lra H_i(G',X;A)$$ defined on the level of chains by $\alpha^*(z)=z$; that is, there is
an inclusion of complexes
$$C_\bullet(X;A)^G\subseteq C_\bullet(X;A)^{G'}.$$  In the particular case of an inclusion of a subgroup $\alpha:H\ra G$, we denote
$\alpha^*$ by
$$\text{res}_H^G:H_i(G,X;A)\lra H_i(H,X;A).$$

\subsection{Homotopy Invariance}\label{homotopy}  Recall that a map $f:X\ra Y$ is a $G$-homotopy equivalence if the map of
fixed-point sets
$$f^H:X^H\lra Y^H$$ is an ordinary homotopy equivalence for every subgroup $H\le G$.  We want to know that the functors $H_i(G,-;A)$
are homotopy invariant.

Our first observation is the following.  A basis of the free abelian group $C_i(X;\zz)^G$ consists of linear combinations of the form
$$\sum_{g\in E_\sigma} g\sigma$$ where $\sigma$ is an $i$-cell of $X$ and $E_\sigma$ is a set of right coset representatives of the
stabilizer $G_\sigma$.  Note that if $G_\sigma = G$, then $\sigma\in C_i(X^G;\zz)$.

\begin{prop} If $f:X\ra Y$ is a $G$-homotopy equivalence, then the chain complexes $C_\bullet(X;\zz)^G$ and $C_\bullet(Y;\zz)^G$ are quasi-isomorphic.
\end{prop}

\begin{proof} Let $h$ be a $G$-homotopy inverse of $f$.  Then $h^H:Y^H\ra X^H$ is a homotopy equivalence for every $H\le G$: $f^Hh^H\simeq\text{id}_{Y^H}$
and $h^Hf^H\simeq\text{id}_{X^H}$.  Consider the induced maps $f_*:C_\bullet(X;\zz)^G\ra C_\bullet(Y;\zz)^G$ and $h_*:C_\bullet(Y;\zz)^G\ra
C_\bullet(X;\zz)^G$.  We show the composite
$$h_*f_*:C_\bullet(X;\zz)^G\lra C_\bullet(X;\zz)^G$$ is chain homotopic to the identity.  If $\sum_{g\in E_\sigma} g\sigma \in C_i(X;\zz)^G$,
then $\sigma\subset X^{G_\sigma}$ and $h^{G_\sigma}f^{G_\sigma}(\sigma) = \sigma + \partial\tau$ for some $\tau\in C_{i+1}(X^{G_\sigma};\zz)$.
It follows that $$h_*f_*\biggl(\sum_{g\in E_\sigma} g\sigma\biggr) = \sum_{g\in E_\sigma}g\sigma + \partial\biggl(\sum_{g\in E_\sigma} g\tau\biggr)$$
and so $h_*f_*\simeq \text{id}$.  The equivalence $f_*h_*\simeq \text{id}$ is proved similarly.
\end{proof}

In the special case $G=\zz/p$, there is an alternate proof that yields more information.  For any finite group $G$ and any
$G$-module $M$, there is a norm map
$$N:M_G\lra M^G$$ defined by $N(\overline{m})=\sum_{g\in G} gm$.  Note that if $M$ is a free $G$-module, then $N$ is an isomorphism.  The norm
map is natural for homomorphisms of $G$-modules, and both $\text{ker}(N)$ and $\text{coker}(N)$ are annihilated by $|G|$.  Consider the norm
map $$N:C_i(X;\zz)_G\lra C_i(X;\zz)^G.$$ Observe that the complex $C_\bullet(X;\zz)_G$ is the cellular chain complex of the quotient space
$X/G$.  Since $C_i(X;\zz)_G$ is a free abelian group, the map $N$ is injective.  If the $G$-action is free, then $N$ is an isomorphism.  We then
have the following.

\begin{prop} If $G$ acts freely on $X$, then the norm map induces an isomorphism
$$H_\bullet(X/G;\zz)\lra H_\bullet(G,X;\zz).$$
\hfill $\qed$
\end{prop}

In the general case, denote the quotient complex $C_\bullet(X;\zz)^G/N(C_\bullet(X;\zz)_G)$ by $D_\bullet(X)$.

\begin{lemma}\label{quotient} If $G=\zz/p$, then $D_\bullet\cong C_\bullet(X^G;\zz/p)$.
\end{lemma}

\begin{proof}  Note that the group $C_k(X;\zz)^G$ has basis consisting of the elements $\sum_{g\in G}g\sigma$
for those $\sigma$ with $G_\sigma=\{1\}$, together with those $\sigma\in C_k(X^G;\zz)$.  The norm map is given by
$$N:\overline{\sigma}\mapsto \sum_{g\in G} g\sigma.$$ If $\sigma\in X^G$, then $\overline{\sigma}=\{\sigma\}$ and $N(\overline{\sigma})
=p\sigma$.  If $G_\sigma=\{1\}$, then $N(\overline{\sigma})=\sum_{g\in G} g\sigma$.  Thus,
$$D_k=\frac{\zz\{\sigma\in (X^G)^{(k)},\sum_{g\in G} g\sigma\}}{\langle \sum_{g\in G} g\sigma, p\sigma\,\text{for}\,\sigma\in X^G\rangle}
\cong C_k(X^G;\zz/p).$$
\end{proof}

\begin{prop}\label{lesprop} If $G=\zz/p$, there is a long exact sequence
$$\ra H_n(X/G;\zz)\ra H_n(G,X;\zz)\ra H_n(X^G;\zz/p)\ra H_{n-1}(X/G;\zz)\ra.$$
\end{prop}

\begin{proof} This is the long exact sequence associated to the exact sequence
$$0\lra C_\bullet(X/G;\zz)\lra C_\bullet(X;\zz)^G\lra D_\bullet(X)\lra 0.$$
\end{proof}

\begin{cor}  If $G=\zz/p$, the functors $H_i(G,-;\zz)$ are homotopy invariant.
\end{cor}

\begin{proof} This is true in general, of course, but in the case of $G=\zz/p$, we
proceed as follows.  If $f:X\ra Y$ is a $G$-homotopy equivalence, then we have homotopy
equivalences
$$\overline{f}:X/G\lra Y/G \qquad \text{and}\qquad f^G:X^G\lra Y^G.$$  The result follows from applying the
Five Lemma to the diagram of long exact sequences obtained via Proposition \ref{lesprop}.
\end{proof}

We shall also make use of the following.  Note that the inclusion of complexes
$$C_\bullet(X;\zz)^G\lra C_\bullet(X)$$
induces a map in homology.  Since an invariant cycle gives rise to an invariant homology class, we really have a map
$$i_*:H_\bullet(G,X;\zz)\lra H_\bullet(X;\zz)^G.$$

\begin{lemma}\label{kernel} The kernel of $i_*$ is annihilated by $|G|$.
\end{lemma}

\begin{proof} Let $[z]\in H_k(G,X;\zz)$ and suppose $i_*([z]) = 0$.  Then $z=\partial\tau$ for some $\tau\in C_{k+1}(X;\zz)$.
Since $z\in C_k(X;\zz)^G$, we have $gz=z$ for every $g\in G$.  Thus, if the $G$-orbit of $\tau$ has length $n$, then
$$\partial\biggl(\sum_{g\in E_\tau} g\tau\biggr) = \sum_{g\in E_\tau}g(\partial\tau) = \sum_{g\in E_\tau} gz = nz.$$
Write $\sigma = \sum_{g\in E_\tau}g\tau$.  Then $g\sigma=\sigma$ for every $g\in G$, and $n[z]=[nz]=[\partial\sigma] = 0$.
Since $n$ divides $|G|$, we have $|G|\cdot\text{ker}(i_*)=0$.
\end{proof}

\subsection{Easy Examples}\label{easy}  As a first example, consider $G=\zz/2$ acting on $X=S^1$ by
$z\mapsto \overline{z}$ (view $S^1$ as the unit circle in ${\mathbb C}$).  We must choose an equivariant
cell decomposition of $X$; this is obtained by taking two $0$-cells ($\{\pm 1\}$) and two $1$-cells
(the upper and lower half-circles).  Denote the vertices by $v_{-1}$ and $v_1$ and the arcs by
$e_+$ and $e_-$.  Then the complex $C_\bullet(X)^G$ is
$$0\lra \zz\{e_+ + e_-\}\stackrel{\partial}{\lra} \zz\{v_{-1},v_1\}\lra 0$$
where $\partial(e_+ + e_-)=2v_{-1} + 2v_1$.  Thus, $H_1(\zz/2,S^1)=0$ and
$$H_0(\zz/2,S^1)=\text{coker}(\partial) \cong \zz\oplus\zz/2.$$ Note that we have the exact sequence
$$0\lra H_1(X/G)\lra H_1(G,X)\lra H_1(X^G;\zz/2) \lra$$
$$\lra H_0(X/G) \lra H_0(G,X)\lra H_0(X^G;\zz/2)\lra 0.$$
Here $X/G \approx [-1,1] \simeq *$ and $X^G=\{v_{-1},v_1\}$.  The map $H_0(X/G)\ra H_0(G,X)$ is
induced by the norm map $C_0(X/G)\ra C_0(X)^G$ and is thus given by $[\overline{v}_1]\mapsto 2[v_1]$.

A similar example is given by $G=\zz/p$ acting on $X=S^1$ by counterclockwise rotation by $2\pi/p$.  Here
the cell decomposition consists of $p$ $0$-cells $v_1,\dots ,v_p$ and $p$ $1$-cells $e_1,\dots ,e_p$.  Since
the $G$-action is free, we expect to recover the homology of $X/G$.  The complex $C_\bullet(X)^G$ is
$$0\lra\zz\{e_1+\cdots +e_p\}\stackrel{\partial}{\lra}\zz\{v_1,\dots ,v_p\}\lra 0$$
with $\partial(e_1+\cdots +e_p)=0$.  So, $H_1(\zz/p,S^1)=\zz$ and $H_0(\zz/p,S^1)=\zz$.  Here, we have
$X/G\approx S^1$ and $X^G=\emptyset$ so that we have the exact sequence
$$0\ra H_1(X/G)\stackrel{\cong}{\ra} H_1(G,X)\ra 0\ra H_0(X/G)\stackrel{\cong}{\ra} H_0(G,X)\ra 0.$$
In this case the maps $H_i(X/G)\ra H_i(G,X)$ are isomorphisms:
$$[\overline{e}_1]\mapsto [e_1+\cdots +e_p] \qquad\text{and}\qquad [\overline{v}_1]\mapsto [v_1+\cdots v_p].$$

A third simple example is given by $G=\zz/2$ acting on $X=S^2$ by $(x,y,z)\mapsto (x,y,-z)$.  An equivariant triangulation
consists of one vertex $v\sim (1,0,0)$, one $1$-cell $e\sim (x,y,0)$ and two $2$-cells $f_+$ and $f_-$ for the northern and
southern hemispheres.  The complex $C_\bullet(X)^G$ is then
$$0\lra\zz\{f_++f_-\}\stackrel{\partial_2}{\lra} \zz\{e\}\stackrel{\partial_1}{\lra}\zz\{v\}\lra 0$$
with $\partial_2(f_++f_-)=2e$ and $\partial_1(e)=0$.  It follows that $H_2(\zz/2,S^2)=0$, $H_1(\zz/2,S^2)\cong\zz/2$,
and $H_0(\zz/2,S^2)\cong \zz$.  Here, $X/G$ is the unit disc and $X^G=S^1$.  The exact sequence is
$$0\ra H_2(\zz/2,S^2)\ra H_2(X^G;\zz/2)\ra H_1(X/G) \ra H_1(\zz/2,S^2)\ra$$
$$\ra H_1(X^G;\zz/2)\ra H_0(X/G)\ra H_0(\zz/2,S^2)\ra H_0(X^G;\zz/2)\ra 0$$
where the maps $H_1(\zz/2,S^2)\ra H_1(X^G;\zz/2)$ and $H_0(X/G)\ra H_0(\zz/2,S^2)$ are the identity and multiplication
by $2$, respectively.

\section{The Functors $S_\bullet(G,X)$ and the Associated Spectral Sequences}\label{auxiliary}

The simple examples in the previous section illustrate what is
going on, but in each case, the quotient space $X/G$ (and its
homology) is easy to describe.  What one really wants is some
means to compute $H_\bullet(G,X)$ without resorting to the chain
complex $C_\bullet(X)^G$, and then to use the map
$H_\bullet(X/G)\ra H_\bullet(G,X)$ to study the homology of $X/G$.
We carry this out in this section.

The fundamental observation is that the invariants functor
$M\mapsto M^G$ is the $0$-th cohomology functor $M\mapsto
H^0(G;M)$. This leads to the following construction.  Let
$P_\bullet\ra\zz$ be a projective resolution over $\zz G$.  Let
$X$ be a finite $G$-CW-complex and consider the fourth quadrant
double complex
$$E^0_{s,t} = \text{Hom}_{\zz G}(P_{-t},C_s(X)).$$
We define groups $S_\bullet(G,X)$ by
$$S_n(G,X)=h_n(\text{Tot}(E^0_{\bullet,\bullet})).$$  Of course, we can do this with any abelian group of coefficients.
We note without proof the following obvious fact.

\begin{prop} The functors $S_\bullet(G,X)$ are independent of the choice of projective resolution $P_\bullet\ra\zz$ or
$G$-equivariant CW-decomposition of $X$. \hfill $\qed$
\end{prop}

As usual, there are two (homological) spectral sequences associated to the double complex.  Since $X$ is assumed to be finite,
convergence is not an issue (there are only finitely many columns).  Taking vertical cohomology, the first spectral sequence has
$${}^{I} E^{1}_{s,t} = H^{-t}(G,C_s(X))\Longrightarrow S_{s-t}(G,X).$$
Note that we obtain
$${}^I E^2_{s,0} = H_s(G,X)$$ so that the groups of interest appear as the $t=0$ row of the ${}^I E^2$-term.

We now analyze the  second spectral sequence. Suppose that
$P_\bullet\ra \zz$ is a projective resolution over $\zz G$. Then
$$E^0_{s,t} = \text{Hom}_{\zz G}(P_{-t},C_s(X)).$$
The spectral sequence is obtained from the row filtration; this
leads to an interchange in the roles of $s$ and $t$ and since the
functors $\text{Hom}_{\zz G}(P_{-t},-)$ are exact, we obtain
$${}^{II} E^1_{s,t} = \text{Hom}_{\zz G}(P_{-s}, H_t(X)) \Longrightarrow S_{t-s}(G,X)$$ with
$d^1:{}^{II}E^1_{s,t}\ra {}^{II}E^1_{s-1,t}$ being the map induced
by the map $P_{-s} \ra P_{-s+1}$ ($s<0$). But this is just the
complex for computing the cohomology of $G$ with coefficients in
$H_t(X)$ and thus
$${}^{II}E^2_{s,t} = H^{-s}(G,H_t(X)).$$ It is useful to visualize
this as living in the second quadrant and the differentials $d^r$ have bidegree
$(-r,r-1)$ as expected.

Since we are primarily interested in the action of $G=\zz/p$, let us restrict attention to this case.  Consider the usual $2$-periodic
resolution
$$\stackrel{N}{\lra}\zz G\stackrel{g-1}{\lra}\zz G\stackrel{N}{\lra}\zz G\stackrel{g-1}{\lra}\zz G\stackrel{\varepsilon}{\lra}\zz$$
where $G=\langle g|g^p=1\rangle$, and $N=\sum_{i=0}^{p-1}g^i$.  The first spectral sequence has
$${}^I E^1_{s,t} = H^{-t}(G,C_s(X))$$ with $d^1$ induced by the boundary map in $C_\bullet(X)$.  These cohomology groups are easily
described.  When $t\ne 0$ is even, we have (see \cite{brown},
p.~58)
$$H^{-t}(G,C_s(X)) = C_s(X)^G/N(C_s(X))$$ and thus, the complex $H^{-t}(G,C_\bullet(X))$ is simply $C_\bullet(X)^G/N(C_\bullet(X))$,
which we denoted by $D_\bullet(X)$ in Section \ref{basic}.  According to Lemma \ref{quotient}, we have $D_\bullet(X)\cong C_\bullet(X^G;\zz/p)$.
Thus, if $t\ne 0$ is even, we have
$${}^I E^2_{s,t} = H_s(X^G;\zz/p).$$
For $t$ odd, we have
\begin{eqnarray*}
H^{-t}(G,C_s(X)) & = & \text{ker}\{\overline{N}:C_s(X)_G\lra C_s(X)^G\} \\
                 & = & 0.
\end{eqnarray*}
However, if we take $\zz/p$ coefficients, then for $t$ even we still have ${}^I E^2_{s,t} = H_s(X^G;\zz/p)$, but for $t$ odd
we obtain
\begin{eqnarray*}
H^{-t}(G,C_s(X;\zz/p)) & = & \text{ker}\{\overline{N}:C_s(X/G;\zz/p)\lra C_s(X;\zz/p)^G\} \\
                       & = & C_s(X^G;\zz/p).
\end{eqnarray*}
Thus, for $t$ odd we also have
$${}^I E^2_{s,t} = H_s(X^G;\zz/p).$$

As an exercise in using these spectral sequences, we prove the following.

\begin{prop} Let $X$ be a finite $G$-CW-complex, where $G$ is any finite group.  Let $A=\zq$ or $A=\zz/\ell$, where $\ell$ does not
divide $|G|$. Then for all $i\ge 0$,
$$H_i(G,X;A) \cong H_i(X;A)^G.$$
\end{prop}

\begin{proof}  Note that we already know this result as the assumption on $A$ allows us to conclude that the norm map
$$C_\bullet(X/G;A)\lra C_\bullet(X;A)^G$$ is an isomorphism, and it is a well-known fact \cite{borel} that with the assumption on $A$,
we have $H_\bullet(X/G;A)\cong H_\bullet(X;A)^G$.  The spectral sequences we have constructed will yield another proof.

Consider the first spectral sequence
$${}^I E^1_{s,t} = H^{-t}(G,C_s(X;A)) \Longrightarrow S_{s-t}(G,X;A).$$
Since $|G|$ is invertible in $A$, the cohomology groups of $G$ vanish in degrees $t\ne 0$.  The spectral sequence therefore
collapses at ${}^I E^2$ and we have
$$S_n(G,X;A) \cong H_n(G,X;A),\, n\ge 0.$$
On the other hand, the second spectral sequence has
$${}^{II} E^2_{s,t} = H^{-s}(G,H_t(X;A)) \Longrightarrow S_{t-s}(G,X;A).$$
Again, these terms vanish for $s\ne 0$ and so the spectral sequence collapses at ${}^{II} E^2$:
$$S_n(G,X;A) \cong H^0(G,H_n(X;A)) = H_n(X;A)^G.$$
\end{proof}

In all the applications that follow, we have $G=\zz/p$.  Since we know how to compute $H_\bullet(G,X;A)$ for $p$ invertible
in $A$, we will primarily use $A=\zz/p$.  The main result we need is the following.

\begin{theorem}\label{collapse} Let $G=\zz/p$ and let $X$ be a finite $G$-CW-complex. Let $A$ be $\zz$ or $\zz/p$.
 Then the spectral sequence ${}^I E^\bullet_{\bullet,\bullet}$,
with $A$ coefficients, satisfies ${}^I E^2 = {}^I E^\infty$.
\end{theorem}

\begin{proof} First note that if $G$ acts trivially on $X$, then
the differential $d^0$ is alternately the zero map and
multiplication by $p$. It follows that with $\zz/p$-coefficients,
all the differentials $d^r$ for $r\ge 2$ vanish. In the general
case, note that the $E^2$-term involves the groups
$H_\bullet(X^G;\zz/p)$ and since $X^G$ is a space on which $G$
acts trivially, the differentials must vanish for $r\ge 2$.

Alternatively, note that the rows in the spectral sequence are
alternately the kernel and cokernel of the norm map, so that when
one lifts elements to compute differentials, one lifts to the zero
element at each stage.
\end{proof}

The interested reader is invited to check that using these
spectral sequences to analyze the examples given in Section
\ref{easy} yields the same calculations.

\section{Smith Theory and the Conner Conjecture}\label{smiththeory}

Smith theory is the generic name given to the study of $G=\zz/p$ actions on homology spheres.  A proof of the following result
using Bredon cohomology may be found in \cite{may}, p.~35.  We use our techniques to provide a different proof.

\begin{theorem}\label{smiththm} Let $X$ be a finite-dimensional $G$-CW-complex, where $G=\zz/p$.  Then
\begin{enumerate}
\item If $H_\bullet(X;\zz/p)\cong H_\bullet(S^n;\zz/p)$ for some $n$, then either $X^G=\emptyset$ or
$H_\bullet(X^G;\zz/p)\cong H_\bullet(S^m;\zz/p)$ for some $m$; and

\item If $H_\bullet(X;\zz/p)\cong H_\bullet(\text{\em pt};\zz/p)$, then $H_\bullet(X^G;\zz/p)\cong
H_\bullet(\text{\em pt};\zz/p)$.
\end{enumerate}
\end{theorem}

\begin{proof}  Consider the first spectral sequence.  By Theorem \ref{collapse}, we have
$${}^I E^\infty_{s,t} = H_s(X^G;\zz/p)$$ for all $t<0$, and ${}^I E^\infty_{s,0}=H_s(G,X;\zz/p)$.  This implies, then,
that for $r<0$ we have
$$S_r(G,X;\zz/p) \cong H_0(X^G;\zz/p)\oplus H_1(X^G;\zz/p)\oplus\cdots\oplus H_m(X^G;\zz/p)$$
(here, $m$ is the largest degree in which $X^G$ has homology; it is certainly finite).  As for the
second spectral sequence, we have
$${}^{II} E^2_{s,t} = H^{-s}(G,H_t(X;\zz/p)).$$  In the case where $X$ is a homology sphere, this consists of two
rows, the $0$-th and the $n$-th, each of which has $H_t(X;\zz/p)=\zz/p$ in each position (note that the action of $G$ on
$H_n(X;\zz/p)$ must be trivial as the only way $\zz/p$ can act on itself is via the identity).  Thus,
$${}^{II} E^2 = {}^{II} E^3 = \cdots {}^{II} E^n$$ and the first possible nontrivial differential is $d^{n+1}$.  Note that
for each $s$, the map $d^{n+1}_{s,0}$ is a map from $\zz/p$ to $\zz/p$ and hence is either $0$ or an isomorphism.  In either
case, we see that for $r\le 0$, $S_r(G,X;\zz/p)$ has rank $0$ or $2$, while for $r>0$, the terms ${}^{II} E^n_{s,n}=H_n(X;\zz/p)$
live to infinity so that $S_r(G,X;\zz/p)\cong \zz/p$ for $0<r\le n$.  Denote by $b_i$ the rank of $H_i(X^G;\zz/p)$.  Then combining
the above we have
$$\sum_{i\ge 0} b_i = 0\;\text{or}\; 2,$$
and hence either $X^G=\emptyset$ or $H_\bullet(X^G;\zz/p)\cong H_\bullet(S^m;\zz/p)$ for some $m$.

In the case where $X$ is acyclic mod $p$, the second spectral sequence has only the $t=0$ row and so $S_r(G,X;\zz/p)\cong
\zz/p$ for $r\le 0$.  It follows then that $\sum b_i = 1$ and hence $X^G$ is acyclic mod $p$ as well.
\end{proof}

We also obtain the following proof of the Conner conjecture, which
was first proved by R.~Oliver \cite{oliver}.

\begin{theorem}\label{conner}  Let $G=\zz/p$ and suppose $X$ is a
finite $G$-CW-complex with $\tilde{H}_\bullet(X;\zz)=0$.  Then
$\tilde{H}_\bullet(X/G;\zz)=0$.
\end{theorem}

\begin{proof}  The second spectral sequence has
$${}^{II}E_{s,t}^2 = H^{-s}(\zz/p,H_t(X;\zz))\Longrightarrow
S_{t-s}(\zz/p,X;\zz).$$  This vanishes for $t>0$ and since $\zz/p$
acts trivially on $H_0(X;\zz)=\zz$, we have
$${}^{II}E_{s,0}^\infty = H^{-s}(\zz/p;\zz)=\begin{cases}
                                              \zz/p & s=2i, i<0 \\
                                               0 &
                                               \text{otherwise}.
                                               \end{cases}$$
Thus, $$S_n(\zz/p,X;\zz) = \begin{cases}
                             \zz & n=0 \\
                             \zz/p & n<0\;\text{even} \\
                              0 & \text{otherwise}.
                              \end{cases}$$

Consider the first spectral sequence:
$${}^IE_{s,t}^\infty = \begin{cases}
                          H_s(X^G;\zz/p) & t\ne 0\;\text{even} \\
                            0 & t\;\text{odd}.
                          \end{cases}$$
Thus, if $n<0$ is even, we have
$$\bop_{i\ge 0} H_{2i}(X^G;\zz/p)\cong S_n(\zz/p,X;\zz)\cong\zz/p,$$
and if $n<0$ is odd, we have
$$\bop_{i\ge 0} H_{2i+1}(X^G;\zz/p)=0.$$
Thus, $\tilde{H}(X^G;\zz/p)=0$ (a fact we could have deduced from
Theorem \ref{smiththm} above), and $H_n(\zz/p,X;\zz)=0$ for $n>0$.
Via the long exact sequence
$$\ra H_n(X/G;\zz)\ra H_n(\zz/p,X;\zz)\ra H_n(X^G;\zz/p)\ra$$
we see that $H_n(X/G;\zz)=0$ for $n>0$, while we have a short
exact sequence
$$0\ra H_0(X/G;\zz)\stackrel{\times p}{\ra} H_0(\zz/p,X;\zz)\ra
H_0(X^G;\zz/p)\ra 0,$$ so that $H_0(X/G;\zz)\cong \zz$.
\end{proof}

\section{Examples}\label{examples}

\subsection{Complex Projective Space}
Our primary example is the following.  Let $G=\zz/2$ act on $X=\cp^n$ by complex conjugation.  We will be interested in the stable
answer ($G$ acting on $\cp^\infty$), but in low homological degrees, we may use a finite projective space to compute the
answer.  This eliminates convergence issues in the spectral sequence calculations.  Note that $X^G=\rp^n$ here.

Consider the first spectral sequence.  By Theorem \ref{collapse}, we have
$${}^I E^\infty_{s,t} = H_s(X^G;\zz/2)$$ for $t<0$ and ${}^I E^\infty_{s,0} = H_s(G,X)$.  Recall that the homology of
$X^G=\rp^n$ with $\zz/2$ coefficients is $\zz/2$ in degrees $0\le i\le n$.  It follows that
$$S_r(\zz/2,\cp^n;\zz/2) \cong (\zz/2)^{n+1},\qquad r<0.$$  Also, for $r\ge 0$, we have
$$S_r(\zz/2,\cp^n;\zz/2) \cong H_r(\zz/2,\cp^n;\zz/2) \oplus (\zz/2)^{n-r}.$$

Now consider the second spectral sequence:
$${}^{II} E^2_{s,t} = H^{-s}(G,H_t(X;\zz/2)).$$
The homology of $\cp^n$ with $\zz/2$ coefficients is $\zz/2$ in even degrees up to $2n$ and $0$ otherwise.  Thus
${}^{II} E^2 = {}^{II} E^3$, and the first possible nontrivial differential is $d^3$.  We claim, however, that
all the differentials vanish.  To see this, note that on the diagonals below the line $t=-s$, we have $n+1$ copies of
$\zz/2$.  As all the differentials involving these terms stay below the main diagonal, and since $S_r(\zz/2,\cp^n;\zz/2)$
has rank $n+1$ for $r<0$, these terms must all live to infinity.  But then all the differentials involving terms on or
above the main diagonal must also vanish.  This follows by noting that the construction of the spectral sequence
involves a certain periodicity; that is, any differential involving a term on or above the main diagonal is the same as some
differential beginning in the same row below the main diagonal.  Therefore, the ${}^{II} E^2$ term is the ${}^{II} E^\infty$
term and we see that for $i\ge 1$
$$S_{2i}(\zz/2,\cp^n;\zz/2) \cong S_{2i-1}(\zz/2,\cp^n;\zz/2) \cong (\zz/2)^{(n+1)-i}.$$

Combining the two calculations above, we see that for $2i\le n$,
\begin{eqnarray*}
H_{2i-1}(\zz/2,\cp^n;\zz/2) & \cong & (\zz/2)^i \\
H_{2i}(\zz/2,\cp^n;\zz/2) & \cong & (\zz/2)^{i+1}
\end{eqnarray*}
while for $2i-1\ge n$,
$$
H_{2i-1}(\zz/2,\cp^n;\zz/2)  \cong H_{2i}(\zz/2,\cp^n;\zz/2) \cong  (\zz/2)^{n+1-i}.$$

Now, if we let $n$ go to infinity, we get the following.

\begin{prop} For all $i\ge 1$,
\begin{eqnarray*}
H_{2i-1}(\zz/2,\cp^\infty;\zz/2) & \cong & (\zz/2)^i \\
H_{2i}(\zz/2,\cp^\infty;\zz/2) & \cong & (\zz/2)^{i+1}.
\end{eqnarray*}
\hfill $\qed$
\end{prop}

We may now use the Universal Coefficient Theorem to calculate the integral groups.  First note
the following fact, which one proves using the ring structure on the cohomology of $\cp^\infty$.

\begin{prop}\label{cpinfty} For all $i\ge 0$,
$$H_i(\cp^\infty;\zz)^{\zz/2} \cong \begin{cases}
                                        0 & i\ne 4k \\
                                        \zz & i=4k.
                                        \end{cases}$$
\hfill $\qed$
\end{prop}

It is also easy to see that $H_{4k}(\cp^\infty;\zz)^{\zz/2}$ is generated by an invariant cycle; that is,
the map
$$i_*:H_{4k}(\zz/2,\cp^\infty;\zz) \lra H_{4k}(\cp^\infty;\zz)^{\zz/2}$$
is surjective.  Moreover, by Lemma \ref{kernel}, we obtain the following.

\begin{prop}\label{kernel2tor} If $i\ne 4k$, then $H_i(\zz/2,\cp^\infty;\zz)$ is $2$-torsion.
If $i=4k$, then $H_i(\zz/2,\cp^\infty;\zz) \cong \zz \oplus A$, where $A$ is a $2$-torsion group. \hfill $\qed$
\end{prop}

Using Proposition \ref{kernel2tor}, in conjunction with the Universal Coefficient Theorem, we obtain the
integral calculation. Details are left to the reader.

\begin{prop}\label{integral} For all $k\ge 1$,
\begin{eqnarray*}
H_{4k-3}(\zz/2,\cp^\infty;\zz) & \cong & (\zz/2)^k \\
H_{4k-2}(\zz/2,\cp^\infty;\zz) & \cong & (\zz/2)^k \\
H_{4k-1}(\zz/2,\cp^\infty;\zz) & \cong & (\zz/2)^k \\
H_{4k}(\zz/2,\cp^\infty;\zz) & \cong & \zz \oplus (\zz/2)^k.
\end{eqnarray*}
\hfill $\qed$
\end{prop}

Of course, $H_0(\zz/2,\cp^\infty;\zz)=\zz$.

Finally, to compute the groups $H_\bullet(\cp^\infty/(\zz/2);\zz)$, we use the exact sequence of Proposition
\ref{lesprop}.  We claim that this splits into short exact sequences
$$0\lra H_i(\cp^\infty/(\zz/2);\zz)\lra H_i(\zz/2,\cp^\infty;\zz)\stackrel{p_*}{\lra} H_i(\rp^\infty;\zz/2)\stackrel{d}{\lra} 0$$
for each $i$, where $d$ is the connecting homomorphism.  Indeed, let $z$ be the generating cycle for $H_i(\rp^\infty;\zz/2)$.
If $i$ is odd, then $z$ is an integral cycle and $z\in C_i(\cp^\infty;\zz)^{\zz/2}$ so that $p_*([z])=[z]$.  If $i$ is even,
then $\partial z=2w$ where $w$ generates $H_{i-1}(\rp^\infty;\zz)$.  But $2[w]=0$ in $H_{i-1}(\cp^\infty/(\zz/2);\zz)$. This follows from
the existence of a transfer map $\text{tr}:H_\bullet(X/G;\zz)\ra H_\bullet(X;\zz)$ for any $G$ and any $G$-space $X$.  In this case, we have that the
composite
$$H_r(\cp^\infty/(\zz/2);\zz) \stackrel{\text{tr}}{\lra} H_r(\cp^\infty;\zz) \lra H_r(\cp^\infty/(\zz/2);\zz)$$
is multiplication by $2$, and if $r$ is odd, we must have $2\cdot H_r(\cp^\infty/(\zz/2);\zz)=0$.
 Thus, $d[z]=[\partial z] = 2[w]=0$, and
so $p_*$ is surjective in this case as well.

Moreover, Proposition \ref{kernel2tor} implies that the above
exact sequences split, and so we obtain the final calculation.

\begin{theorem}  For all $k\ge 1$,
\begin{eqnarray*}
H_{4k-3}(\cp^\infty/(\zz/2);\zz) & \cong & (\zz/2)^{k-1} \\
H_{4k-2}(\cp^\infty/(\zz/2);\zz) & \cong & (\zz/2)^{k-1} \\
H_{4k-1}(\cp^\infty/(\zz/2);\zz) & \cong & (\zz/2)^{k-1} \\
H_{4k}(\cp^\infty/(\zz/2);\zz) & \cong & \zz \oplus (\zz/2)^{k-1}.
\end{eqnarray*}
\hfill $\qed$
\end{theorem}

Note that this agrees with the calculation of $H^\bullet(\cp^\infty/(\zz/2);\zz)$ given by Dugger in \cite{dugger}.

\subsection{Dihedral Groups}  In \cite{knudson}, the author considered the case of a finite group $G$ acting on the classifying
space $BQ$ of some group $Q$.  For example, if $\zz/2$ acts on $\zz/n$ by $x\mapsto -x$, then one obtains the following
computation for $n\not\equiv 0\mod 4$:
$$H_\bullet(\zz/2,B(\zz/n);\zz) \cong H_\bullet(\zz/n;\zz)^{\zz/2}.$$
The case $n\equiv 0\mod 4$ is more difficult, however.  One first reduces to the calculation in the case of $n=2^s$, $s\ge 2$.  With
mod $2$ coefficients, the Poincar\'e series of $H_\bullet(\zz/2,B(\zz/2^s);\zz/2)$ is
$$p(t) = \frac{1}{(1-t)^2}.$$
What is interesting is that this is also the Poincar\'e series of the mod 2 homology of the dihedral group $D_{s+1}$ of order
$2^{s+1}$:
$$0\lra \zz/2^s \lra D_{s+1} \lra \zz/2\lra 0$$
(for a proof, see \cite{mitchell}).  Thus, the complex $C_\bullet(B(\zz/2^s);\zz/2)^{\zz/2}$ ``computes" the homology of the whole
dihedral group.  Of course, the composition
$$H_\bullet(\zz/2,B(\zz/2^s);\zz/2)\lra H_\bullet(\zz/2^s;\zz/2)\lra H_\bullet(D_{s+1};\zz/2)$$
is not an isomorphism as the first map has a large kernel (the middle term is $\zz/2$ in all degrees, while the first
term is much larger).  This is a rather curious phenomenon.

\end{document}